\documentclass[12pt,draft]{amsart}
\usepackage{amsmath,amsthm,latexsym,amscd,amsbsy,amssymb,fleqn,leqno}
\chardef\bslash=`\\ % p. 424, TeXbook

\makeatletter
\def\verbatim{\interlinepenalty\@M \@verbatim
  \leftskip\@totalleftmargin\advance\leftskip2pc
  \frenchspacing\@vobeyspaces \@xverbatim}
\makeatother
\makeatletter
  \def\dgt@k{\dg@DX=-3 \dg@DY=2 \dg@SIZE=3} 
\makeatother

\makeatletter
  \def\dgt@kk{\dg@DX=3 \dg@DY=-1 \dg@SIZE=3}%
\makeatother

\theoremstyle{plain}
\newtheorem{thm}{Theorem}[section]
\newtheorem{cor}[thm]{Corollary}
\newtheorem{lem}[thm]{Lemma}
\newtheorem{pro}[thm]{Proposition}

\theoremstyle{definition}

\numberwithin{equation}{section}

\newcounter{rmnum}

\def\symbolnote#1#2{\let\thefootn=\thefootnote%
\renewcommand{\thefootnote}{\fnsymbol{footnote}}%
\footnotemark[#1]%
\footnotetext[#1]{#2}%
\let\thefootnote=\thefootn
}

\newfont{\bbb}{msbm10 scaled \magstep1}
\newfont{\bbc}{msbm8 scaled \magstep0}

\newcommand{\R}{\mbox{\bbb R}}

\newcommand{\N}{\mbox{\bbb N}}

\newcommand{\uin}{\mbox{\bbb I}}

\begin{document}

%%%%%%% Begin Topmatter %%%%%%%%%%

\title[On finite-dimensional  maps]{On finite-dimensional maps}
\author{H. Murat Tuncali}
\address{Department of Mathematics,
Nipissing University\-wan,
100 College Drive, P.O. Box 5002, North Bay, ON, P1B 8L7, Canada}
\email{muratt@nipissingu.ca}
\thanks{The first author was partially supported by NSERC grant.}

\author{Vesko Valov}
\address{Department of Mathematics, Nipissing University,
100 College Drive, P.O. Box 5002, North Bay, ON, P1B 8L7, Canada}
\email{veskov@nipissingu.ca}
\thanks{The second author was partially supported by Nipissing University Research Council Grant.}

 
%%%%%%% End topmatter %%%%%%%%%

\begin{abstract}{Let $f\colon X\to Y$ be a perfect  surjective map of metrizable spaces. It is shown that if $Y$ is  a $C$-space
(resp., $\dim Y\leq n$ and $\dim f\leq m$), then the function space $C(X,\uin^{\infty})$ (resp.,  $C(X,\uin^{2n+1+m})$) equipped with the source limitation topology contains a dense $G_{\delta}$-set $\mathcal{H}$ such that $f\times g$ embeds $X$ into $Y\times\uin^{\infty}$ (resp., into $Y\times\uin^{2n+1+m}$)  for every $g\in\mathcal{H}$.  Some applications of this result are also given.}
\end{abstract}

\maketitle
\markboth{H. M.~Tuncali and V.~Valov}{On finite-dimensional maps}

%%%%%%%%%%%%%%%%%%%%%%%%%%%%%%%%%%%%%%%%%%%%%%%%%%%%%%%%%%%%

\section{Introduction}

The aim of this note is to prove the following theorem.

\begin{thm}
Let $f\colon X\to Y$ be a perfect surjection of countable functional weight between paracompact spaces.  
\begin{itemize}
\item[(i)] If  $\dim f\leq n$ and $\dim Y\leq m$, then the set $\mathcal{H}$ of all maps $g\in C(X,\uin^{2n+1+m})$ 
with $f\times g$ being an embedding is dense and $G_{\delta}$ in $C(X,\uin^{2n+1+m})$ with respect to the source limitation topology.
\item[(ii)] If $Y$ is a $C$-space, then all maps $g\colon X\to\uin^{\infty}$ such that $f\times g$ is an embedding form a dense and $G_{\delta}$-subset of $C(X,\uin^{\infty})$ with respect to the source limitation topology.
\end{itemize}
\end{thm}

Here, $f$ is of countable weight \cite{bp:98} if there exists a map $h\colon X\to\uin^{\infty}$ such that $f\times g$ embeds $X$ into
$Y\times\uin^{\infty}$. The $C$-space property was introduced by Haver \cite{wh:74} for compact metric spaces and then extended by Addis and Gresham \cite{ag:78} for general spaces. The source limitation topology (or the fine topology) is defined in Section 2, it is well known that this topology is stronger than the uniform convergence topology.  We consider the covering dimension $\dim$ (see \cite{re:95}) and $\dim f\leq n$ means that $\dim f^{-1}(y)\leq n$ for every $y\in Y$.

Theorem 1.1(ii) was announced without proof in \cite{tv:01}.  Concerning (i), 
Eilenberg \cite{se:35} proved that $\mathcal{H}\neq\emptyset$ for $n=0$ and compact metric spaces $X$ and $Y$. 
Then Pasynkov \cite{bp:67} extended Eilenberg's result to perfect $0$-dimensional maps between arbitrary metric spaces and proved the compact version of Theorem 1.1(i) in \cite{bp:96}.
The Eilenberg theorem was also generalized by S. Bogatyi, V. Fedorchuk and J. van Mill \cite{bfm:00}.  Recently, Pasynkov \cite{bp:98} extended his results from \cite{bp:96}. It follows from \cite[Theorem 9.5]{bp:98} that, under the hypotheses of Theorem 1.1(i), the set $\mathcal{H}$ is dense in $C(X,\uin^{2n+1+m})$ with respect to the uniform convergence topology generated by the Euclidean  metric on $\uin^{2n+1+m}$.  Since the function space $C(X,\uin^{2n+1+m})$  with any reasonable topology (e.g., the uniform convergence and the fine topology) has Baire property, in many cases it is important to know that $\mathcal{H}$ is also residual (i.e., contains a dense $G_{\delta}$-set) in $C(X,\uin^{2n+1+m})$. We demonstrate that importance by proving the following generalization of a result due to Roberts \cite{r:41} (see Proposition 4.3):
\begin{cor}
Under the hypothesis of Theorem $1.1(i)$, there exists a dense $G_{\delta}$-set $\mathcal{B}$ in  
$C(X,\uin^{2n+1+m})$ with the source limitation topology such that $\mathcal{B}\subset\mathcal{H}$ and every $g\in\mathcal{B}$ has the following property:
$\dim g(f^{-1}(y))\cap\Pi\leq k-n-1$ for any $y\in Y$  and any rational $k$-plane $\Pi$ in $\R^{2n+1+m}$ with $ 0\leq k\leq 2n$.
\end{cor}

Another application of Theorem 1.1(i) is an improvement of the mentioned above result of  Bogatyi, V. Fedorchuk and  van Mill \cite{bfm:00} (see Proposition 4.1).  Let us also note that Theorem 1.1(i) extends a result of Hansen \cite{h:78} 
about embedding of finite covers.   

The paper is organized as follows.  Section 2 is devoted to the main technical result, Theorem 2.1. Theorem 1.1 is proved in Section 3 and the final Section 4 contains some applications.  All single-valued maps under discussion are continuous, and all function spaces, if not explicitely stated otherwise, are equipped with the source limitation topology.

%%%%%%%%%%%%%%%%%%%%%%%%%%%%%%%%%
\section{$\omega$-maps}

For any spaces $X$ and $Y$ by $C(X,Y)$ we denote the set of all continuous maps from
$X$ into $Y$. If $(Y,d)$ is a metric space, then the source limitation topology on $C(X,Y)$ is defined in the following way:
a subset $U\subset C(X,Y)$ is open in $C(X,Y)$ with respect to
the source limitation topology provided for every $g\in U$ there exists 
a continuous function $\alpha\colon X\to (0,\infty)$ such that $\overline{B}(g,\alpha)\subset U$. Here, $\overline{B}(g,\alpha)$ denotes the set  
$\{h\in C(X,Y):d(g(x),h(x))\leq\alpha (x)\hbox{}~~\mbox{for each 
$x\in X$}\}$. 

The source limitation topology is also known as the fine topology and $C(X,Y)$ with this topology has Baire property provided $(Y,d)$ is a complete metric space \cite{jm:75}.  Moreover, the source limitation topology on $C(X,Y)$ doesn't depend on the metric of $Y$ when $X$ is paracompact \cite{nk:69}.
 
Throughout the paper $\uin^k$ denotes the $k$-dimensional cube equipped with the Euclidean metric $d_k$, and $D_k$ denotes the uniform convergence metric on $C(X,\uin^k)$ generated by $d_k$. 
We also agree to denote by $cov(X)$ the family of all open covers of $X$. In case $(X,d)$ is a metric space,
$B_{\epsilon}(x)$ (resp., $\overline{B}_{\epsilon}(x)$) stands for the
open (resp., closed) ball in $(X,d)$ with center $x$ and radius
$\epsilon$.

It is well known that if $X$ is a compact $n$-dimensional space, then for every $\omega\in cov(X)$ the set of all $\omega$-maps $g\colon X\to\uin^{2n+1}$ is open and dense in $C(X,\uin^{2n+1})$ (recall that $g\in C(X,\uin^{2n+1})$ is an $\omega$-map if  for every $x\in X$ there is a neighborhood $V$ of $g(x)$ in  $\uin^{2n+1}$ such that 
$g^{-1}(V)\subset U$ for some $U\in\omega$).  Next theorem, which easily implies Theorem 1.1(i),  is a counterpart of the above result for $n$-dimensional maps.   

\begin{thm}
Let $f\colon X\to Y$ be a perfect surjection between paracompact spaces such that $\dim f\leq n$ and  $\dim Y\leq m$. Then, 
for every $\omega\in cov(X)$, the set $\{g\in C(X,\uin^{2n+1+m}): f\times g \hbox{ is an $\omega$-map}\}$ is 
open and dense in $C(X,\uin^{2n+1+m})$.
\end{thm}

The proof of this theorem follows very closely the proof of Theorem 2.2 from \cite{tv:01}. We need few lemmas, in all these lemmas we suppose that $f$, $X$, $Y$ are as in Theorem 2.1, $\omega\in cov(X)$ and $k=2n+1+m$. 
We also denote by $C(X,Y\times\uin^k,f)$
the set of all maps $h\colon X\to Y\times\uin^k$ such that $\pi_Y\circ h=f$, where
$\pi_Y\colon Y\times\uin^k\to Y$ is the projection. For any closed $K\subset X$, $C_{\omega}(X|K,Y\times\uin^k,f)$ stands for the set of all $h\in C(X,Y\times\uin^k,f)$
with $h|K$ being an $\omega$-map (as a map from $K$ into $Y\times\uin^k$) and 
$C_{\omega}(X|K,\uin^k)$ consists of all $g\in C(X,\uin^k)$ such that 
$f\times g\in C_{\omega}(X|K,Y\times\uin^k,f)$. In case $K=X$ we simply write 
$C_{\omega}(X,Y\times\uin^k,f)$ (resp., $C_{\omega}(X,\uin^k)$) instead of
$C_{\omega}(X|X,Y\times\uin^k,f)$ (resp., $C_{\omega}(X|X,\uin^k)$). 

\begin{lem}
Let $g\in C_{\omega}(X|f^{-1}(y),\uin^k)$ for some $y\in Y$. Then there exists a neighborhood $U$ of $y$ in $Y$ such that the restriction $g|f^{-1}(U)$ is an $\omega$-map. Moreover,  $g\in C_{\omega}(X,\uin^k)$ provided $g\in C_{\omega}(X|f^{-1}(y),\uin^k)$ for every $y\in Y$.
\end{lem}

\begin{proof} Since $f$ is perfect, all fibers of $f$ are compact and the map $f\times g\colon X\to Y\times\uin^k$ is closed.  These two facts imply the proof (see Lemma 2.3 and Corollary 2.4 from \cite{tv:01} for similar arguments).
\end{proof}    

\begin{lem}
For any closed $K\subset X$  the set 
$C_{\omega}(X|K,\uin^k)$ is open in $C(X,\uin^k)$. 
\end{lem}

\begin{proof}
Let $g_0\in C_{\omega}(X|K,\uin^k)$. We are going to find $\alpha\in C(X,(0,\infty))$ with $\overline{B}(g_0,\alpha)\subset C_{\omega}(X|K,\uin^k)$. For every $y\in H=f(K)$ there exists a neighborhood $U_y$  of $y$ in $Y$ such that $g_0|(f^{-1}(U_y)\cap K)$ is an $\omega$-map. Then $\omega_1=\{U_y:y\in H\}\cup\{Y\backslash H\}$ is an open cover of $Y$. Using that $Y$ is paracompact, we can find a metric space $(M,d)$, a surjection $p\colon Y\to M$ and $\mu\in cov(M)$ such that $p^{-1}(\mu)$ refines $\omega_1$.
Hence, every $z\in p(H)$ has a neighborhood $W_z$ in $M$ such that $g_0|(p\circ f)^{-1}(W_z)\cap K$ is an $\omega$-map. The last condition implies that $h_0|K$ also is an $\omega$-map, where 
$h_0=(p\circ f)\times g_0$. The proof of next claim is similar to the proof of the Claim from \cite[Lemma 2.5]{tv:01}.

\medskip
{\em Claim.} {\em There exists a locally finite open family $\gamma$ in $M\times\uin^k$ covering $h_0(K)$ such that
every $g\in C(X,\uin^k)$ belongs to $C_{\omega}(X|K,\uin^k)$ provided $h|K$ is $\gamma$-close to $h_0|K$, where $h=(p\circ f)\times g$.}

\medskip 

Let $\rho$ be the metric on $M\times\uin^k$ defined by $\rho (t_1,t_2)=d(z_1,z_2)+d_k(w_1,w_2)$, where $t_i=(z_i,w_i)$, $i=1,2$. Denote by $V$ the union of all elements of the family $\gamma$ from the claim let $\alpha_1 \colon K\to (0,\infty)$ be the continuous function 
$\alpha_1 (x)=2^{-1}\sup\{\rho (h_0(x),V\backslash W):W\in\gamma\}$. If 
$h=(p\circ f)\times g$ with $g\in C(X,\uin^k)$ and
$\rho(h_0(x),h(x))\leq\alpha_1(x)$ for every $x\in K$, then $h|K$ is $\gamma$-close to 
$h_0|K$. According to the claim, the last relation yields that 
$g\in C_{\omega}(X|K,\uin^k)$. We take a continuous extension $\alpha\colon X\to (0,\infty)$ of $\alpha_1$ and observe that $d_k(g_0(x),g(x))=\rho(h_0(x),h(x))$ for every $x\in X$. 
Therefore, $\overline{B}(g_0,\alpha)\subset C_{\omega}(X|K,\uin^k)$.
\end{proof}

\begin{lem}
If $C(X,\uin^k)$ is equipped with the uniform convergence topology, then the set-valued map $\displaystyle\psi_{\omega}\colon Y\to 2^{C(X,\uin^k)}$, defined by the formula
$\psi_{\omega}(y)=C(X,\uin^k)\backslash C_{\omega}(X|f^{-1}(y),\uin^k)$, has a closed graph.
\end{lem}

\begin{proof} See the proof of \cite[Lemma 2.6]{tv:01}
\end{proof}

\begin{lem}
Let $K$ be a compact space such that $\dim K\leq n$. 
Then for every $\gamma\in cov(K)$ the set of all maps $g\in C(\uin^m\times K,\uin^k)$ with each $g|(\{z\}\times K)$, $z\in\uin^n$, being a
$\gamma$-map $($as a map from $K$ into $\uin^k$$)$ 
is dense in $C(\uin^m\times K,\uin^k)$ with respect to the uniform convergence topology.
\end{lem}

\begin{proof}
We fix $\gamma\in cov(K)$, a map $g_0\colon\uin^m\times K\to\uin^k$ and $\epsilon>0$. If $K$ is metrizable, the Pasynkov version of Theorem 1.1(i) (for compact spaces) \cite{bp:96} yields a map $g\colon \uin^m\times K\to\uin^k$ which is $\epsilon$-close to   $g_0$ and $\pi\times g$ embeds $\uin^m\times K$ into $\uin^m\times\uin^k$, where $\pi\colon \uin^m\times K\to\uin^m$ denotes the projection  onto $\uin^m$.  Then, each $g|(\{z\}\times K)$,  $z\in\uin^n$, is a $\gamma$-map because it embeds
$ \{z\}\times K$ into $\uin^k$.  

If $K$ is not metrizable we can 
represent $K$ as the limit space of a $\sigma$-complete inverse system $\mathcal S=\{K_{\beta},p_{\beta}^{\beta +1}:\beta\in B\}$ such that each $K_{\beta}$ is a metrizable compactum with $\dim K_{\beta}\leq n$. 
Applying standard inverse spectra arguments (see \cite{book}), we can find $\theta\in B$, $\gamma _1\in cov(K_{\theta})$ and $g_{\theta}\in C(\uin^m\times K_{\theta},\uin^k)$ such that $g_{\theta}\circ (id\times p_{\theta})=g_0$ and $p_{\theta}^{-1}(\gamma _1)$ refines $\gamma$, where $p_{\theta}\colon K\to K_{\theta}$ denotes the $\theta$th limit projection and $id$  the identity map on $\uin^m$. Then, by virtue of the previous case, there exists a map 
$g_1\in C(\uin^m\times K_{\theta},\uin^k)$ which is $\epsilon$-close to $g_{\theta}$ and the restrictions 
$g_1|(\{z\}\times K_{\theta})$ are $\gamma_1$-maps. It follows from our construction 
that $g=g_1\circ (id\times p_{\theta})$ is $\epsilon$-close to $g_0$ and each
$g|(\{z\}\times K)$ is an $\gamma$-map.
\end{proof}

Recall that a closed subset $F$ of the metrizable apace $M$ is said to be a $Z$-set in $M$ \cite{vm:89},   if the set $C(\uin^{\infty},M\backslash F)$ is dense in $C(\uin^{\infty},M)$ with respect to the uniform convergence topology. If, in the above definition, $\uin^{\infty}$ is replaced by $\uin^m$, $m\in\N\cup\{0\}$, we say that
$F$ is a $Z_m$-set in $M$.

\begin{lem}
Let $\alpha\colon X\to (0,\infty)$ be a positive continuous function and $g_0\in C(X,\uin^k)$.
Then $\psi_{\omega}(y)\cap\overline{B}(g_0,\alpha)$ is a $Z_m$-set in $\overline{B}(g_0,\alpha)$ for every $y\in Y$, where $\overline{B}(g_0,\alpha)$ is considered as a subspace of $C(X,\uin^k)$ with the uniform convergence topology. 
\end{lem}

\begin{proof}
We follow the proof of \cite[Lemma 2.8]{tv:01}. 
All function spaces in this  proof are equipped with the uniform convergence topology.
By Lemma 2.4, each $\psi_{\omega}(y)$
is closed in $C(X,\uin^k)$. Hence, $\psi_{\omega}(y)\cap\overline{B}(g_0,\alpha)$ is closed in $\overline{B}(g_0,\alpha)$. We need to show that, for fixed $y\in Y$, $\delta>0$ and a map $u\colon\uin^m \to \overline{B}(g_0,\alpha)$ there exists a map
$v\colon\uin^m\to\overline{B}(g_0,\alpha)\backslash\psi_{\omega}(y)$ which is $\delta$-close to $u$ with respect to the uniform metric $D_k$. To this end, observe first that $u$ generates $h\in C(\uin^m\times X,\uin^k)$, $h(z,x)=u(z)(x)$, such that
$d_k(h(z,x),g_0(x))\leq\alpha (x)$ for any $(z,x)\in\uin^m\times X$. Since $f^{-1}(y)$ is compact, there exists $\lambda\in (0,1)$ such that $\lambda\sup\{\alpha (x):x\in f^{-1}(y)\}<\displaystyle\frac{\delta}{2}$. Define $h_1\in C(\uin^m\times f^{-1}(y),\uin^k)$ by $h_1(z,x)=(1-\lambda)h(z,x)+\lambda g_0(x)$. Then, for every $(z,x)\in\uin^m\times f^{-1}(y)$, we have \\

\smallskip\noindent   
(1) \hbox{}~~~~~~$d_k(h_1(z,x),g_0(x))\leq (1-\lambda)\alpha (x)<\alpha (x)$ \\

\smallskip\noindent
and

\smallskip\noindent
(2) \hbox{}~~~~~~$d_k(h_1(z,x),h(z,x))\leq\lambda\alpha (x)<\displaystyle\frac{\delta}{2}$.

\smallskip\noindent
Let $\displaystyle q<\min\{r,\frac{\delta}{2}\}$, where $r$ is the positive number
$\inf\{\alpha (x)-d_k(h_1(z,x),g_0(x)):(z,x)\in\uin^m\times f^{-1}(y)\}$. 
Since $\dim f^{-1}(y)\leq n$, by Lemma 2.5 (applied to the product $\uin^m\times f^{-1}(y)$),   
there is a map $h_2\in C(\uin^m\times f^{-1}(y),\uin^k)$ such that $d_k(h_2(z,x),h_1(z,x))<q$ and $h_2|(\{z\}\times f^{-1}(y))$ is an $\omega$-map for each $(z,x)\in\uin^m\times f^{-1}(y)$. Then, by $(1)$ and $(2)$, for all $(z,x)\in\uin^m\times f^{-1}(y)$ we have \\

\smallskip\noindent
(3) \hbox{}~~~~~~$d_k(h_2(z,x),h(z,x))<\delta$ and $d_k(h_2(z,x),g_0(x))<\alpha (x)$. \\

\smallskip\noindent
Because both $\uin^m$ and $f^{-1}(y)$ are compact, $u_2(z)(x)=h_2(z,x)$ defines 
the map $u_2\colon\uin^m\to C(f^{-1}(y),\uin^k)$. The required map $v$ will be obtained as a lifting of $u_2$.
Let $\pi\colon\overline{B}(g_0,\alpha)\to C(f^{-1}(y),\uin^k)$ denote the restriction map $\pi (g)=g|f^{-1}(y)$. 
As in the proof of \cite[Lemma 2.8]{tv:01}, we can see that
$u_2(z)\in\pi (\overline{B}(g_0,\alpha))$ for every $z\in\uin^m$ and
$\theta(z)=\overline{\pi^{-1}(u_2(z))\cap B_{\delta}(u(z))}$ defines a convex-valued map from $\uin^m$ into $\overline{B}(g_0,\alpha)$ which is lower semi-continuous. By  the Michael selection theorem \cite[Theorem 3.2"]{em:56}, $\theta$ has a continuous selection $v\colon\uin^m\to C(X,\uin^k)$. Then $v$ maps $\uin^m$ into $\overline{B}(g_0,\alpha)$ and $v$ is $\delta$-close to $u$. Moreover, for any $z\in\uin^m$ we have $\pi(v(z))=u_2(z)$ and $u_2(z)$, being the restriction $h_2|(\{z\}\times f^{-1}(y))$, is an $\omega$-map. Hence, 
$v(z)\in C_{\omega}(X|f^{-1}(y),\uin^k)$, $z\in\uin^m$, i.e.
$v\colon\uin^m\to\overline{B}(g_0,\alpha)\backslash\psi_{\omega}(y)$.       
\end{proof}

\begin{lem}
The set $C_{\omega}(X,\uin^k)$ is dense in $C(X,\uin^k)$.
\end{lem}

\begin{proof}
It suffices to show that, for fixed $g_0\in C(X,\uin^k)$ and a positive continuous function $\alpha\colon X\to (0,\infty)$, there exists $g\in \overline{B}(g_0,\alpha)\cap C_{\omega}(X,\uin^k)$. The space $C(X,\uin^k)$ with the uniform convergence topology 
is a closed convex subspace of the the Banach space $E$ consisting of all bounded continuous maps from $X$ into $\R^k$.
We define the set-valued map $\phi$ from $Y$ into $C(X,\uin^k)$, 
$\phi(y)=\overline{B}(g_0,\alpha)$, $y\in Y$. According to Lemma 2.6, 
$\overline{B}(g_0,\alpha)\cap\psi_{\omega}(y)$ is a $Z_m$-set in $\overline{B}(g_0,\alpha)$ for every $y\in Y$. So, we have a lower semi-continuous closed and convex-valued map $\phi\colon Y\to 2^E$ and another map $\psi_{\omega}\colon Y\to 2^E$ with a closed graph (see Lemma 2.4) such that 
$\phi(y)\cap\psi_{\omega}(y)$ is a $Z_m$-set in $\phi(y)$ for each $y\in Y$. Moreover, $\dim Y\leq n$, so we can apply 
\cite[Theorem 1.2]{gv:01} to obtain a  continuous map $h\colon Y\to C(X,\uin^k)$ with $h(y)\in\phi(y)\backslash\psi_{\omega}(y)$ for every $y\in Y$ (we actually apply the following partial case of \cite[Theorem 1.2]{gv:01}: If $\Phi\colon Z\to 2^F$ is a lower semi-continuous closed and convex-valued map
from a paracompact space $Z$ with $\dim Z\leq n$ into a Banach space $E$  and $\psi\colon Z\to 2^E$ is a set-valued map with a closed graph such that $\Phi(z)\cap\psi(z)$ is a $Z_m$-set in $\Phi(z)$ for all $z\in Z$, then there exists a map $q\colon Z\to E$ with
$q(z)\in\Phi(z)\backslash\psi(z)$ for each $z\in Z$).  
Observe that $h$ is a map from $Y$ into $\overline{B}(g_0,\alpha)$ such that $h(y)\not\in\psi_{\omega}(y)$ for every $y\in Y$, i.e. 
$h(y)\in \overline{B}(g_0,\alpha)\cap C_{\omega}(X|f^{-1}(y),\uin^k)$, $y\in Y$. Then
$g(x)=h(f(x))(x)$, $x\in X$, defines a map $g\in \overline{B}(g_0,\alpha)$ such that  $g|f^{-1}(y)=h(y)|f^{-1}(y)$ for all $y\in Y$. Hence,
$g\in C_{\omega}(X|f^{-1}(y),\uin^k)$, $y\in Y$. Finally, by virtue of Lemma 2.2, $g\in C_{\omega}(X,\uin^k)$.  
\end{proof} 

%%%%%%%%%%%%%%%%%%%%%%%%%%%%%%%%%

%%%%%%%%%%%%%%%%%%%%%%%%%%%%%%
%%%%%%%%%%%%%%%%%%%%%%%%%%%%%%

\section{Proof of Theorem 1.1}

We first prove Theorem 1.1(i). Fix a map $h\colon X\to\uin^{\infty}$ such that $f\times h$ embeds $X$ into $Y\times\uin^{\infty}$, and let $\omega_i=h^{-1}(\gamma_i)$, where $\{\gamma_i\}$ is a sequence of open covers of $\uin^{\infty}$ with $\displaystyle\hbox{mesh} (\gamma_i)<2^{-i}$ for every $i$.   
Then, by Theorem 2.1, each $\displaystyle C_{\omega_i}(X,\uin^{2n+1+m})$ is open and dense in $C(X,\uin^{2n+1+m})$, so 
$\displaystyle C_0=\bigcap_{i=1}^{\infty}C_{\omega_i}(X,\uin^{2n+1+m})$ is dense and $G_{\delta}$ in $C(X,\uin^{2n+1+m})$.  According to the choice of $\omega_i$, $ C_0\subset\mathcal{H}$.  
On the other hand, any $g\in\mathcal{H}$ belongs to $C_0$. Therefore, $C_0=\mathcal{H}$ and we are done.

The proof of Theorem 1.1(ii) follows the same scheme as that one of (i). The only difference is that, instead of Theorem 2.1, we use the following proposition.
\begin{pro}
Let $f\colon X\to Y$ be a perfect surjection such that $Y$ is a paracompact $C$-space . Then, 
for every $\omega\in cov(X)$, the set of all $g\in C(X,\uin^{\infty})$ with $f\times g$ being an $\omega$-map is 
open and dense in $C(X,\uin^{\infty})$.
\end{pro}
To prove Proposition 3.1, we fix an $\omega\in cov(X)$, an embedding of $\uin^{\infty}$ in $\displaystyle l_2$ as a closed convex subset and let $d_{\infty}$ denote the metric on $\uin^{\infty}$ generated by the norm of $\displaystyle l_2$.  Following the proof of Theorem 2.1, we show how to modify the hypotheses and the proofs of Lemmas 2.2 - 2.7.  One general modification is to replace everywhere the metric $d_k$ by $d_{\infty}$ and both $\uin^m$ and $\uin^k$ by $\uin^{\infty}$. With this general correction, the
modified versions of Lemmas 2.2 -
2.4 remain valid. The modification of Lemma 2.5 should read as follows: Let $K$ be a compact space. Then, for every $\gamma\in cov(K)$, the set of all maps $g\in C(\uin^{\infty}\times K,\uin^{\infty})$ with each $g|(\{z\}\times K)$, $z\in\uin^{\infty}$, being a $\gamma$-map is dense in 
$C(\uin^{\infty}\times K,\uin^{\infty})$.  In the proof of this new version of Lemma 2.5 we need the following corrections (except the general one): to produce the map $g$ in the case when $K$ is metrizable, choose $g$ to be an embedding of $\uin^{\infty}\times K$ into $\uin^{\infty}$ which is $\epsilon$-close to $g_0$;  when $K$ is not metrizable, we omit the restriction $K_{\beta}$ to be at most $n$-dimensional.  In Lemma 2.6 (hypothesis and the proof) we need to replace $Z_m$-set by $Z$-set. 
The new version of Lemma 2.7 is: The set $C_{\omega}(X,\uin^{\infty})$ is dense in $C(X,\uin^{\infty})$.  Here is a sketch of the modified proof: For given 
$\alpha\colon X\to (0,\infty)$ and $g_0\in C(X,\uin^{\infty})$ we define $\phi\colon X\to 2^E$ to be the constant map $\phi (y)=\overline{B}(g_0,\alpha)$, where $E$ is now the Banach space of all bounded maps from $X$ into $\displaystyle l_2$. Then, according to the modified versions of Lemmas 2.4 and 2.6, $\psi_{\omega}$ has a closed graph and $\phi(y)\cap\psi_{\omega}(y)$ is a $Z$-set in $\phi(y)$ for every $y\in Y$.  Next, apply \cite[Theorem 1.1]{gv:99} to obtain a map $h\colon Y\to C(X,\uin^{\infty})$ with $h(y)\in\phi(y)\backslash\psi_{\omega}(y)$, $y\in Y$. Finally, the required map $g\in C_{\omega}(X,\uin^{\infty})$ is defined by $g(x)=h(f(x)(x)$.

%%%%%%%%%%%%%%%%%%%%%%%%%%%%%%%%%

%%%%%%%%%%%%%%%%%%%%%%%%%%%%%%
%%%%%%%%%%%%%%%%%%%%%%%%%%%%%%

\section{Some applications}

Let $f\colon X\to Y$ be a map and $B\subset X$ an $F_{\sigma}$-subset of $X$. We say that the restriction $f|B$ is $\sigma$-perfect with a countable weight if  $B$ is the union of closed subsets $B_i\subset X$, $i\in\N$, such that all restrictions $f|B_i\colon B_i\to f(B_i)$ are perfect maps with countable weight and $f(B_i)\subset Y$ closed. 
Our first application is the following generalization of \cite[Theorem 4 and Theorem 5]{bfm:00} (following the notations of Proposition 4.1 below, it was shown in \cite{bfm:00} that the sets $\mathcal{A}$ and $\mathcal{A}_{\infty}$ are non-empty in the more restrictive situation when both $X$ and $Y$ are separable metric spaces and $B$ is $\sigma$-compact).  

\begin{pro}
Let $f\colon X\to Y$ be a surjective map between paracompact spaces and $B\subset X$ an $F_{\sigma}$-set. Suppose $f|B$ is $\sigma$-perfect with countable weight such that $\dim (f|B)\leq n$ and $\dim f(B)\leq m$. Then we have:
\begin{itemize}
\item[(i)]
The set 
$\mathcal{A}$ consisting of all $g\in C(X,\uin^{2n+1+m})$ with $(f\times g)|B$ being one-to-one is dense and $G_{\delta}$ in  $C(X,\uin^{2n+1+m})$.
\item[(ii)]
The set $\mathcal{A}_{\infty}$ of all $g=(g_1,g_2,..,)\in C(X,\uin^{\infty})$ with $\displaystyle f\times g_{i_1}\times....\times g_{i_{2n+1+m}}$  being one-to-one on $B$ for any $2n+1+m$ distinct integers $\displaystyle i_1,...,i_{2n+1+m}$ is dense and $G_{\delta}$ in $C(X,\uin^{\infty})$.
\end{itemize}
\end{pro}

\begin{proof}
Let $B=\cup_{i=1}^{\infty}B_i$ such that, for each $i$, $B_i\subset X$, $f(B_i)\subset Y$ are closed sets and $f_i=f|B_i$ is a perfect map with countable weight. Since $\dim (f_i)\leq n$, $\dim f(B_i)\leq m$ and both $B_i$ and $f(B_i)$ are paracompact, by Theorem 1.1(i), every set  $\mathcal{H}_i=\{g\in C(B_i,\uin^{2n+1+m}): f_i\times g \hbox{ is an embedding}\}$ is dense and $G_{\delta}$ in $C(B_i,\uin^{2n+1+m})$.  The restriction maps $p_i\colon C(X,\uin^{2n+1+m})\to C(B_i,\uin^{2n+1+m})$, $p_i(g)=g|B_i$, are open and continuous surjections. Hence, all sets $p^{-1}(\mathcal{H}_i)$ are dense and $G_{\delta}$ in $C(X,\uin^{2n+1+m})$, so is $\mathcal{A}=\cap_{i=1}^{\infty}p^{-1}(\mathcal{H}_i)$.  

To prove (ii), for any set $\Gamma$ of $2n+1+m$ distinct integers $\displaystyle i_1,...,i_{2n+1+m}$ we represent $\uin^{\infty}$ as the product
$\uin^{\Gamma}\times\uin^{N\backslash\Gamma}$.  It is easy to show that 
$C(X,\uin^{\infty})$ is homeomorphic to $C(X,\uin^{\Gamma})\times C(X,\uin^{N\backslash\Gamma} )$.  Let $\pi_{\Gamma}$ denote the projection from $C(X,\uin^{\infty})$ onto $C(X,\uin^{\Gamma})$ and 
$\mathcal{A}_{\Gamma}=\{g\in C(X,\uin^{\Gamma}): (f\times g)|B \hbox{ is one-to-one } \}$. Then, by (i), each 
$\pi_{\Gamma}^{-1}(\mathcal{A}_{\Gamma})$ is dense and $G_{\delta}$ in $C(X,\uin^{\infty})$. Since $\mathcal{A}_{\infty}$ is the intersection of all   $\pi_{\Gamma}^{-1}(\mathcal{A}_{\Gamma})$, we are done.  
\end{proof}

\begin{lem}
Let $H=\cup_{i=1}^{\infty}H_i$ be an $n$-dimensional subset of the paracompact space $X$ with each $H_i\subset X$ closed.  Suppose $E$ is a Banach space and $Y\subset E$ a closed convex subset. Then, for every $\sigma Z_n$-set $K$ in $Y$, the set of all $g\in C(X,Y)$ with $g(H)\cap K=\emptyset$ is dense and $G_{\delta}$ in $C(X,Y)$. 
\end{lem}

\begin{proof}
Since the restriction maps from $C(X,Y)$ into $C(H_i,Y)$ are continuous and open, it suffices to show that, for any $i$, 
the set $\mathcal{G}_i$ of all $g\in C(H_i,Y)$ with $g(H_i)\cap K=\emptyset$ is dense and $G_{\delta}$ in $C(H_i,Y)$. Let 
$K=\cup_{j=1}^{\infty}K_j$, where each $K_j$ is a closed $Z_n$-subset of $Y$, and $\mathcal{G}_{ij}$ consist of all
$g\in C(H_i,Y)$ with $g(H_i)\cap K_j=\emptyset$.  Then $\mathcal{G}_i=\cap_{j=1}^{\infty}\mathcal{G}_{ij}$.  Therefore, our proof is further reduced to the proof that $\mathcal{G}_{ij}$ is open and dense in $C(H_i,Y)$ for any $i$ and $j$.  Finally, observe that, for fixed $i, j\in\N$, we can assume $H_i=X$ and $K_j=K$.  Under this assumption, we are going now to show that the set 
$U=\{g\in C(X,Y): g(X)\cap K=\emptyset\}$ is open and dense in $C(X,Y)$.

Suppose $g_0\in U$ and let $\alpha (x)=2^{-1}d(g_0(x),K)$, where $d$ is the metric on $Y$ generated by the norm of $E$. Then 
$\overline{B}(g_0,\alpha)\subset U$, so $U$ is open in $C(X,Y)$. To show that $U$ is dense, take arbitrary $h\in C(X,Y)$ and $\alpha\in C(X,(0,\infty))$ and consider the set-valued map $\phi\colon  X\to 2^Y$, $\phi (x)=\overline{B}_{\alpha (x)}(h(x))$. Then $\phi$ is lower semi-continuous with closed and convex values. Moreover, since $K$ is a $Z_n$-set in $Y$, $K\cap B_{\alpha (x)}(h(x))$ is a $Z_n$-set in $B_{\alpha(x)}(h(x))$, $x\in X$. This implies that 
$\phi (x)\cap K$ is a $Z_n$-set in $\phi (x)$ for every $x\in X$ (see, for example, \cite[Lemma 2.3]{gv:99}). Therefore, we can apply \cite[Theorem 1.2]{gv:01} to obtain a map $g\colon X\to Y$ such that $g(x)\in\phi (x)\backslash K$ for all $x$. Hence, $g\in U\cap\overline{B}(g,\alpha)$ which shows that $U$ is dense in $C(X,Y)$. 
\end{proof}

Let $\displaystyle N^{2n+1+m}_{n+m}$ denote the set of all points from $\uin^{2n+1+m}$ having at most $n+m$ rational coordinates. It is well known that the compliment of $\displaystyle N^{2n+1+m}_{n+m}$ in $\uin^{2n+1+m}$ is the union of countably many $n$-dimensional subspaces $M_i$ with each $M_i$ being a $Z_{n+m}$-set in $\uin^{2n+1+m}$. Moreover, if $X$ is a paracompact space of dimension $\dim X\leq n+m$, according to Lemma 4.2, each of the sets 
$\mathcal{M}_i=\{g\in C(X,\uin^{2n+1+m}): g(X)\cap M_i=\emptyset\}$ is open and dense in $C(X,\uin^{2n+1+m})$. This implies the following generalization of the N\"{o}beling-Pontryagin-Lefschetz theorem:  In the hypothesis of Theorem 1.1(i), the set $\mathcal{N}$ of all
$g\in\mathcal{H}$ such that  $g(X)\subset N^{2n+1+m}_{m+n}$ is
dense and $G_{\delta}$ in $C(X,\uin^{2n+1+m})$. Indeed, by the generalized Hurewicz theorem \cite{es:62}, $\dim X\leq n+m$, so the sets $\mathcal{M}_i$ defined above are open and dense in $C(X,\uin^{2n+1+m})$.  Then $\mathcal{N}$, being the intersection of all $\mathcal{M}_i$  and $\mathcal{H}$, is dense and $G_{\delta}$ in $C(X,\uin^{2n+1+m})$.   But we can prove a stronger result, which generalizes a theorem due to J. Roberts \cite{r:41} as well. Below, by a (rational) $k$-plane in $\R^p$ we mean a subspace of the form $\displaystyle\{(x_1,..,x_p)\in\R^p: x_{i_j}=r_j, j=1,..,p-k\}$, where $\{i_1,..,i_{p-k}\}$ is a subset of $\{1,..,p\}$ and $r_j$, $j=1,..,p-k$, are fixed (rational) numbers. Obviously, every $k$-plane is a $k$-dimensional affine subspace of $\R^p$. 

\begin{pro}
Under the hypothesis of Theorem $1.1(i)$, $\mathcal{H}$ contains a dense $G_{\delta}$-subset $\mathcal{B}$ of 
$C(X,\uin^{2n+1+m})$ such that every $g\in\mathcal{B}$ satisfies the following condition $($R$)$:
$\dim g(f^{-1}(y))\cap\Pi\leq k-n-1$ for any $y\in Y$  and any rational $k$-plane $\Pi$ in $\R^{2n+1+m}$ with $ 0\leq k\leq 2n$.
\end{pro}

\begin{proof}
By \cite{bp:98}, there exists a map $h\colon X\to\uin^n$ such that $f\times h$ is $0$-dimensional.  Take 
 a cover $\{D_i: i=1,..,n+1\}$ of $\uin^n$ consisting of $0$-dimensional $G_{\delta}$-sets and let $A_i=(f\times h)^{-1}(Y\times D_i)$.  Obviously, each $A_i$ is $G_{\delta}$ in $X$ and $X=\cup_{i=1}^{n+1}A_i$. Moreover, according to the choice of $h$, all restrictions $f|A_i$ are $0$-dimensional.  Now, for any $n+1\leq k\leq 2n$ consider the set 
$\displaystyle B_k=X\backslash (A_1\cup...\cup A_{k-n})$.  Then 
$\displaystyle B_k\subset A_{k-n+1}\cup...\cup A_{n+1}$ which implies $\dim f^{-1}(y)\cap B_k\leq 2n-k$ for every $y\in Y$ (recall that  each $f|A_i$ is $0$-dimensional). On other hand, $B_k$ is $F_{\sigma}$ in $X$, so is $f(B_k)$ in $Y$. Consequently, $\dim f(B_k)\leq m$. Applying the generalized Hurewicz theorem to the $\sigma$-perfect map $f|B_k$, we conclude that $\dim B_k\leq m+2n-k$.  If $0\leq k\leq n$, we denote $B_k=X$. Then $\dim B_k\leq m+n\leq m+2n-k$ (by the generalized Hurewicz theorem). Hence, for any $k=0,1,..,2n$ we have $\dim B_k\leq m+2n-k$.

Next, let $\Pi (k)$ be the union of all rational $k$-planes in $\R^{2n+1+m}$. Since any $k$-plane $\Pi$ is a $Z_{m+2n-k}$-set in $\R^{2n+1+m}$, so is $\Pi\cap\uin^{2n+1+m}$ in $\uin^{2n+1+m}$.   
Hence, $\Pi (k)\cap\uin^{2n+1+m}$ is a
$\sigma Z_{m+2n-k}$-set in $\uin^{2n+1+m}$. Then, by Lemma 4.2, each of the sets 
$\mathcal{B}_k=\{g\in C(X,,\uin^{2n+1+m}): g(B_k)\cap\Pi (k)=\emptyset\}$, $k=0,1,...,2n$, is dense and $G_{\delta}$ in $C(X,\uin^{2n+1+m})$. 
Therefore, the intersection $\mathcal{B}$ of $\mathcal{H}$ and all $\mathcal{B}_k$, $0\leq k\leq 2n$, is 
dense and $G_{\delta}$ in $C(X,\uin^{2n+1+m})$.
It remains only to show that every $g\in\mathcal{B}$ satisfies the condition (R).
This holds if $0\leq k\leq n$ because $g\in\mathcal{B}_k$, so $g(X)$ doesn't meet any rational $k$-plane. If $n+1\leq k\leq 2n$, then $g\in\mathcal{H}$ implies that $g$ embeds every fiber $f^{-1}(y)$, $y\in Y$, in $\uin^{2n+1+m}$. On the other hand, 
$g(B_k)\cap\Pi (k)=\emptyset$, so $g(f^{-1}(y))\cap\Pi (k)$ is homeomorphic to a subset of
$f^{-1}(y)\backslash B_k$. Since $f^{-1}(y)\backslash B_k$ is contained in 
$f^{-1}(y)\cap (A_1\cup...\cup A_{k-n})$ and the last set is a separable metric space of dimension $\leq k-n-1$ (recall that each $f^{-1}(y)\cap A_i$ is $0$-dimensional), we have
$\dim g(f^{-1}(y))\cap\Pi (k)\leq k-n-1$. Consequently, $\dim g(f^{-1}(y))\cap\Pi\leq k-n-1$ for any rational $k$-plane $\Pi$.
\end{proof}

We now consider the case when the map $f$ from Theorem 1.1 is closed but not necessary perfect.

\begin{pro}
Let $f\colon X\to Y$ be a closed $n$-dimensional surjection of countable weight such that $X$ is normal and $Y$ a paracompact  $m$-dimensional space.  Then all maps $g\in C(X,\uin^{2n+1+m})$ such that  $f\times g$ is an embedding and $g$ satisfies condition $($R$)$ from Proposition $4.3$ form a dense subset in  $C(X,\uin^{2n+1+m})$ with respect to the uniform convergence topology.
\end{pro}

\begin{proof}
We fix a map $h\colon X\to\uin^{\infty}$ such that $f\times h$ is an embedding (recall that $f$ is of countable weight) and a sequence $\{\gamma_i\}$ of open covers of $\uin^{\infty}$  with $\displaystyle\hbox{mesh} (\gamma_i)<2^{-i}$. 
We also consider  the \v{C}ech-Stone extension $\beta f$ of $f$ and the space
$Z=(\beta f)^{-1}(Y)$.  Then $Z$ is paracompact (as a perfect preimage of $Y$) and, since  
$X$ is normal,  
 $\bar{f}\colon Z\to Y$ is a perfect $n$-dimensional map. Using this fact, as in the proof of Proposition 4.3, we cover $Z$ by $G_{\delta}$-sets $A_1,..,A_{n+1}$ with each $\bar{f}|A_i$ being $0$-dimensional, and construct the sets $B_k$ for every $0\leq k\leq 2n$. Then, all
sets $\mathcal{B}_k=\{g\in C(Z,\uin^{2n+1+m}): g(B_k)\cap\Pi (k)=\emptyset\}$, $k=0,1,...,2n$, are dense and $G_{\delta}$ in $C(Z,\uin^{2n+1+m})$. Next,
consider the extension $\bar{h}\colon Z\to\uin^{\infty}$ of $h$ and the covers $\omega_i=\bar{h}^{-1}(\gamma_i)\in cov(Z)$. By Theorem 2.1, $\displaystyle C_{\omega_i}(Z,\uin^{2n+1+m})$ are dense and $G_{\delta}$-sets in $C(Z,\uin^{2n+1+m})$, so is the intersection $\mathcal{B}$ of all $C_{\omega_i}(Z,\uin^{2n+1+m})$ and $\mathcal{B}_k$. Then, the set $\mathcal{Q}=\{g|X:g\in\mathcal{B}\}$ may not be dense and $G_{\delta}$ in $C(X,\uin^{2n+1+m})$, but it is dense in $C(X,\uin^{2n+1+m})$ with respect to the uniform convergence topology generated by the Euclidean metric on $\uin^{2n+1+m}$. Moreover, it follows from the choice of the coves $\omega_i$ that 
every $f\times q$, $q\in\mathcal{Q}$, is an embedding. Finally, as in Proposition 4.3, one can show that all $q\in\mathcal{Q}$ satisfy condition (R). 
\end{proof}

\bigskip

\end{document}